\documentclass[11pt]{amsart}

\usepackage{amsfonts,amsmath}
\usepackage{amssymb, latexsym}
\usepackage{amscd,amsthm}
\input xy
\xyoption{all}

\newcommand{\marginlabel}[1]%
  {\mbox{}\marginpar{\raggedleft\hspace{0pt}\bfseries\sf#1}}

\def\CC{{\mathbb C}}
\def\PP{{\textbf P}}
\def\OO{{\mathcal O}}
\def\R{\mathbf{R}}

\def\D{\mathbf{D}}

\def\F{\mathcal{F}}
\def\E{\mathcal{E}}

\def\G{\mathcal{G}}

\def\Pic0{{\rm Pic}^0(X)}

\theoremstyle{plain}

\newtheorem{theorem}{Theorem}[section]
\newtheorem{theoremalpha}{Theorem}

\newtheorem{proposition/example}[theorem]{Proposition/Example}
\newtheorem{proposition/definition}[theorem]{Proposition/Definition}
\newtheorem{proposition}[theorem]{Proposition}
\newtheorem{corollary}[theorem]{Corollary}
\newtheorem{lemma}[theorem]{Lemma}

\theoremstyle{definition}
\newtheorem{definition}[theorem]{Definition}
\newtheorem{remark}[theorem]{Remark}

\newtheorem{conjecture/question}[theorem]{Conjecture/Question}

\newtheorem{remark/definition}[theorem]{Remark/Definition}
\newtheorem{definition/notation}[theorem]{Definition/Notation}

 \theoremstyle{remark}

\begin{document}

\title{Strong Generic Vanishing and a higher dimensional Castelnuovo-de Franchis inequality}

\author[G. Pareschi]{Giuseppe Pareschi}
\address{Dipartamento di Matematica, Universit\`a di Roma, Tor Vergata, V.le della
Ricerca Scientifica, I-00133 Roma, Italy} \email{{\tt
pareschi@mat.uniroma2.it}}

\author[M. Popa]{Mihnea Popa}
\address{Department of Mathematics, University of Illinois at Chicago,
851 S. Morgan Street, Chicago, IL 60607, USA } \email{{\tt
mpopa@math.uic.edu}}

\thanks{MP was partially supported by the NSF grant DMS 0601252
and by a Sloan Fellowship.}

\maketitle

\setlength{\parskip}{.1 in}

\markboth{G. PARESCHI and M. POPA} {\bf A higher dimensional Castelnuovo-de Franchis inequality}

\section{Introduction}

A classical theorem of Castelnuovo \cite{cast} and de Franchis \cite{defranchis}, 
proved at the turn of the last century, implies that a smooth compact complex
surface $X$ satisfying the inequality
\begin{equation}\label{classical} 
p_g(X)< 2q(X)-3
\end{equation}
admits a fibration (i.e. a surjective morphism with connected fibres) onto
a smooth curve of genus at least $2$ (cf. also the modern age references \cite{bpv},
Proposition 4.2 and \cite{beauville1}, C.8). This basic result in surface theory
concerns in fact surfaces whose Albanese map is generically finite:
otherwise the Albanese map of $X$ automatically provides a fibration
onto a smooth curve of genus at least $2$, unless $q(X)\le 1$ when $2q(X)-3<0$.  
Inequality (\ref{classical}) can be restated as
$$\chi(\omega_X)<q(X)-2.$$
The purpose of this paper is to show that, in this form, the inequality of Castelnuovo and de Franchis
has the following natural generalization in arbitrary dimension:

\begin{theoremalpha}\label{cdf} 
Let $X$ be a compact K\"ahler manifold of maximal Albanese dimension. If
$$\chi(\omega_X)<q(X)-\dim X$$
then $X$ admits a surjective morphism with connected fibers onto a normal compact
analytic variety $Y$ with $0 < \dim Y<\dim X$ and with any smooth model $\tilde Y$ of Albanese general type.
Moreover, one can also impose that every $\tilde Y$ be in addition of general type, and with $\chi(\omega_{\tilde Y})> 0$.
\end{theoremalpha}

We recall that a compact K\"ahler manifold is said to be of \emph{maximal Albanese dimension} 
if its Albanese map is generically finite onto its image. According to Catanese \cite{catanese}, it is said to be of \emph{Albanese general type} if in addition its Albanese map is not surjective. Moreover, a morphism to a variety $Y$ as in the theorem is called in \cite{catanese} a \emph{higher irrational pencil}. Thus Theorem \ref{cdf} can be restated as 
follows: \emph{if $X$ is of maximal Albanese dimension and does not admit a higher irrational pencil, then $\chi(\omega_X) \ge q(X) - \dim X$.}  An extremely interesting fact proved in \cite{catanese} Theorem 2.25, inspired by results of Beauville and Siu for fibrations onto curves, is that the non-existence of higher irrational 
pencils is a purely topological condition detected by the cohomology algebra $H^*(X, \CC)$.
Note also that irregular manifolds not of maximal Albanese dimension by definition admit morphisms to lower-dimensional 
nontrivial varieties (with smooth models) of maximal Albanese dimension. It is easy however to construct examples
of such varieties of dimension $3$ or higher satisfying $\chi(\omega_X) < q(X) - \dim X$  and which 
do not admit higher irrational pencils (cf. \S5), so the maximal Albanese dimension hypothesis cannot be dropped in 
Theorem \ref{cdf}.

Curves of genus at least $2$, which appear in the classical result on surfaces, can be generalized
in higher dimension (up to taking smooth models) in at least three different meaningful ways:
(a) varieties of Albanese general type (as already noted in \cite{catanese});   
(b) varieties of maximal Albanese dimension 
whose holomorphic Euler-characteristic is positive; (c) varieties of general type.  
As suggested in the statement of Theorem \ref{cdf}, satisfying the key condition (a) means in fact that all three conditions can be satisfied at the same time, by invoking 
structure results of Ueno and Ein-Lazarsfeld. Theorem \ref{cdf} is in fact a corollary of the following more general result:

\begin{theoremalpha}\label{general} 
Let $X$ be a compact K\"ahler manifold of maximal Albanese dimension. Then $X$ admits a surjective morphism with connected fibers onto a normal compact analytic variety $Y$ with $0\le \dim Y<\dim X$, any smooth model $\tilde Y$ of maximal Albanese dimension, and such that
$$\chi(\omega_X)\ge q(X)- \dim X - (q(\tilde Y)- \dim \tilde Y).$$ 
\end{theoremalpha}

Here we allow the case of maps to a point, i.e. $\dim Y=0$. Except for the last statement in parenthesis, Theorem \ref{cdf} follows immediately from this. Indeed, the hypothesis of Theorem \ref{cdf} implies that 
 $q(\tilde Y)-\dim  \tilde Y > 0$. In particular $\dim \tilde Y>0$ and $\tilde Y$ is of Albanese general type. Note that  in fact we have the more precise lower bound
$$q(\tilde Y)\ge \dim \tilde Y +(q(X)-\dim X -\chi(\omega_X)).$$ 
It only remains to check the final assertion, which is standard: assume that $\tilde Y$ is not of general
type. Then, by Ueno's theorem (\cite{mori} (3.7)), the image $g(Y)$ is ruled by
positive dimensional subtori of $\widehat{T}$ over a basis $Z$ which is a subvariety of
general type of an abelian variety. Consider the Stein factorization $X \rightarrow Z^{\prime}$
of the induced morphism $X \rightarrow Z$. By a result of Ein-Lazarsfeld (\cite{el}, Theorem
3) we have that $\chi (\omega_{\tilde Z^{\prime}}) > 0$ for any smooth model as well, since otherwise $Z$ would have to be further ruled by non-trivial subtori, in contradiction with the fact that it is of general type. Thus we can replace $Y$ by $Z^\prime$, which satisfies all the required conditions. Finally, 
Theorem \ref{general} provides in fact a refinement of the result of Ein-Lazarsfeld \cite{el} mentioned above, 
namely Corollary \ref{ruling} below.

The proof of Theorem \ref{general} relies on an algebraic framework based on results in homological commutative algebra, which holds in great generality, combined with the theory of Generic Vanishing for the canonical bundle (cf. \cite{gl1}, \cite{gl2}, \cite{hacon}). We adopt an approach often employed in 
recent literature on irregular varieties, namely the use of the integral transform for coherent sheaves 
$\R\Phi_P(\F) :=\R {p_{\Pic0}}_*(p_X^* \F\otimes P)$, where $P$ is a Poincar\'e line bundle on $X\times \Pic0$ 
(cf. e.g. \cite{hp3}, \cite{pp1} or \cite{pp2} for a small sampling). We begin by giving in \S2 a 
positive answer, in the compact K\"ahler case, to the natural generalization of a problem 
of Green-Lazarsfeld, \cite{gl2}, Problem 6.2. In the projective case, this was proved by Hacon \cite{hacon} (see also  \cite{pareschi} for a different proof), and extended further in \cite{pp3}.

\begin{theoremalpha}\label{vanishing}
Let $X$ be a compact K\"ahler manifold of dimension $d$, with Albanese image $a(X)$ of dimension 
$d-k$. Then $R^i \Phi_P (\OO_X) (= R^i {p_{\Pic0}}_* (P)) = 0$ for all $i \not\in [d-k, d]$.
\end{theoremalpha}

The proof given here simply observes that the result is equivalent to the Generic Vanishing theorem of \cite{gl1}. This equivalence is a particular case of a result which holds for an arbitrary coherent sheaf on $X$ (cf. Theorem \ref{general_equiv}).
We then introduce the \emph{Generic Vanishing Index}, a measure of the 
size of the cohomological support loci of a sheaf $\F$ on an  irregular variety $X$. 
This index  turns out to quantify a homological property of 
$\R \Phi_P (\F)$, namely that of being a $k$-th syzygy sheaf (cf. Corollary \ref{syzygy_kahler} and 
Definition \ref{gv_index}). Based on the Evans-Griffith Syzygy Theorem, 
a deep result in local commutative algebra, it provides as a consequence a lower bound for the holomorphic Euler characteristic (cf. Theorem \ref{chi_inequality}).\footnote{
We note that all of these results hold in the general context of arbitrary integral functors defined by locally free kernels, characterizing the filtration of ${\rm Coh}(X)$ by $GV_k$-sheaves. The definitions and statements are provided at the end of \S5 for completeness.}
Therefore ``strong" Generic Vanishing, meaning high Generic Vanishing index or equivalently small 
cohomological support loci, makes the Euler characteristic large.
Finally, this implies Theorem \ref{general} in combination with the geometric description of the cohomological support loci of the canonical bundle of a compact K\"ahler manifold given by Green and Lazarsfeld \cite{gl2}.

It may be of interest to note that this approach in arbitrary dimension does not rely on
generalizations of the well-known Castelnuovo-de Franchis Lemma on holomorphic $1$-forms, which 
is the key ingredient in the classical proof on surfaces.  Various highly interesting higher-dimensional 
generalizations of this result do exist in the literature (cf. \cite{catanese}, \cite{gl2}, \cite{ran}).
Their main thrust is however in a different direction, not directly linked to the holomorphic 
Euler characteristic. Finaly, some instances of the inequality we prove here have already been discovered in the nice work of Hacon-Pardini on the classification of certain irregular varieties. The case $\chi (\omega_X) = 1$ is worked out in \cite{hp3}, while some of its consequences already appeared in \cite{hp2}, both relying on a homological approach in \cite{el} with its roots in ideas of Green.

\noindent 
{\bf Acknowledgements.} We thank Lawrence Ein, Christopher Hacon, Mart\'i Lahoz, Gian Pietro
Pirola and Claire Voisin for valuable discussions. Special thanks are due to the referees, whose comments 
have helped us improve the exposition.

\section{Vanishing of higher direct images and Generic Vanishing}

Let $X$ be a smooth projective variety over an algebraically closed field, or a compact complex manifold, with Albanese map $a: X \rightarrow {\rm Alb}(X)$. Consider a Poincar\'e
line bundle $P$ on $X \times {\rm Pic}^0 (X)$. We will make use of the integral functor $\R\Phi_P$ as in the Introduction,
and also the analogous functor $\R\Phi_{P^\vee}$ defined by the dual of $P$. Given a sheaf $\F$ on $X$, we consider the \emph{cohomological support loci}
$$V^i(\F)=\{[\alpha]\in \Pic0\>|\>h^i(\F\otimes\alpha^{-1})>0\}.$$
(Here we denote by $[\alpha]$  a point in $\Pic0$ and by $\alpha=P_{|X\times[\alpha]}$ the corresponding line bundle on $X$.) We consider also the loci $${\rm Supp}~R^i\Phi_{P^\vee}(\F) \subseteq \Pic0.$$ 
Although, by  base-change, ${\rm Supp}~R^i\Phi_{P^\vee}(\F)$ is only contained in $V^i(\F)$, and in
general not equal to it, the two carry the same numerical information (around any point in $\Pic0$) in the following sense:

\begin{lemma}\label{codimension}
For any integer $m$, the following are equivalent for any $[\alpha]\in \Pic0$:

\noindent
(i) ${\rm codim}_{[\alpha]}~{\rm Supp}~R^i \Phi_{P^\vee} (\F) \geq i+m$ for all $i>0$. 

\noindent
(ii) ${\rm codim}_{[\alpha]}~ V^i (\F)\geq i+m$ for all $i>0$.
\end{lemma}
\begin{proof}
All the arguments are local, so we omit mentioning this in what follows.
Since by base-change we have that ${\rm Supp}~R^i \Phi_{P^\vee}
(\F)\subseteq V^i (\F)$, it is enough to prove that (i) implies (ii).
Assume, by contradiction, that (i) holds but (ii) fails. We can find a largest integer $i$ such that there is an irreducible component $W$ of $V^i(\F)$ with ${\rm codim}~ V^{i+1} (\F)< i + m$. By the maximality of $i$, we conclude that $H^{i+1}(\F\otimes \alpha^{-1})=0$ for a general point $[\alpha]\in W$. By  base-change,\footnote{We are using the following well-known base-change result (applied to our setting): if $h^{i+1}(\F\otimes \alpha^{-1})$ is constant in a neighborhood of $X$, then both 
$R^{i+1} \Phi_{P^\vee} (\F)$ and $R^i \Phi_{P^\vee} (\F)$ have the base-change property in a neighborhood of $[\alpha]$. By semicontinuity, this holds if $h^{i+1}(\F\otimes\alpha^{-1})=0$ (see e.g. \cite{mumford}, Cor.2 p.52). } $R^i \Phi_{P^\vee} (\F)$ is non-zero at a general point of $W$, and this contradicts (i). 
\end{proof}

Theorem \ref{vanishing} in the Introduction follows from the more general local result below, as the 
main result of \cite{gl1}, Theorem 1,  
says that for a compact K\"ahler manifold $X$ of dimension $d$ and Albanese dimension $d-k$ we have ${\rm codim} ~V^i (\omega_X) \ge i - k$ for all $i$. We use the notation $\R\Delta \F := \R {\mathcal H}om_X(\F,\omega_X)$. 

\begin{theorem}\label{general_equiv}
Let $X$ be a smooth projective variety over an algebraically closed field, or a compact complex manifold, of dimension $d$, 
and let $\F$ be a coherent sheaf on $X$. The following are equivalent, for any $[\alpha]\in\Pic0$:

\noindent
(i) ${\rm codim}_{[\alpha]} V^i (\F) \ge i - k, {\rm ~for~ all~} i > 0$.

\noindent
(ii) $R^i \Phi_P(\R\Delta \F)_{[\alpha]} = 0$  for all $i \not\in [d-k, d]$.
\end{theorem}
\begin{proof}
By Lemma \ref{codimension}, the condition in (i) is equivalent to
$${\rm codim}_{[\alpha]}~{\rm Supp}~R^i\Phi_{P^\vee}(\F)\ge i -k$$
for all $i>0$. Now Grothendieck duality exchanges the duality and integral functors according to the 
formula
\begin{equation}\label{gd}
\R \mathcal{H}om_{\Pic0} (\R \Phi_{P^\vee} (\F), \OO_{\Pic0}) \cong \R \Phi_P (\R\Delta \F) [d], 
\end{equation}
where $[d]$ denotes a shift to the left by $d$.
(In the context of smooth varieties, this is well known, cf. for example \cite{pp3} Lemma 2.2.
The same proof works  on complex manifolds, due to the fact that the analogue of Grothendieck 
duality holds in that context also, by \cite{rr} and \cite{rrv}).

\noindent
To prove that (i) implies (ii), we make use of the spectral sequence 
$$E_2^{pq} : = \mathcal{E}xt^p _{\Pic0}(R^q\Phi_{P^\vee}(\F), \OO_{\Pic0}) \Rightarrow \mathcal{E}xt^{p - q} _{\Pic0}(\R \Phi_{P^\vee} (\F), \OO_{\Pic0}).$$
In other words, to converge to $\mathcal{E}xt^{i}_{\Pic0}(\R \Phi_{P^\vee} (\F), \OO_{\Pic0})$, we have to start with the $E_2$-terms given by $\mathcal{E}xt^{i +j} _{\Pic0}(R^j \Phi_{P^\vee} (\F), \OO_{\Pic0})$. Since  ${\rm codim}_{[\alpha]}~ {\rm Supp} ~R^j \Phi_{P^\vee} (\F) \ge j-k$ for all $j$, 
by Lemma \ref{ext} we get
$\mathcal{E}xt^{i +j}_{\Pic0} (R^j \Phi_{P^\vee} (\F), \OO_{\Pic0})_{[\alpha]} = 0$ for all $i < -k$, which by the spectral sequence implies that 
$\mathcal{E}xt^{i} _{\Pic0}(\R \Phi_{P^\vee} (\F), \OO_{\Pic0})_{[\alpha]} = 0$ for $i < -k$. By (\ref{gd}), this is the same as saying that $R^i\Phi_P(\R\Delta \F)_{[\alpha]} = 0$ for all $i < d-k$. On the other hand, by Serre-Grothendieck duality, 
$$H^i(\R\Delta \F\otimes \alpha)\cong H^{d-i}(\F\otimes\alpha^{-1})^\vee$$ 
which is $0$ if $i>d$, and so $R^i\Phi_P(\R\Delta \F) = 0$ for all $i>d$ by base-change.\footnote{We recall that
base-change works more generally for hypercohomology of bounded complexes of sheaves which are flat over the base
(\cite{ega3} 7.7, especially 7.7.4, and Remarque 7.7.12(ii))}

\noindent
To prove that (ii) implies (i),\footnote{This type of argument was already used in 
\cite{hacon}, and later in \cite{pp3}.} note that we have an involution 
$$\R \mathcal{H}om_{\Pic0} \big(\R \mathcal{H}om _{\Pic0}(\R \Phi_{P^\vee} (\F), \OO_{\Pic0}), \OO_{\Pic0}\big) \cong \R \Phi_{P^\vee} (\F),$$
which induces a spectral sequence 
$$E_2^{pq} : = \mathcal{E}xt^p_{\Pic0} \big(\mathcal{E}xt^q _{\Pic0}(\R \Phi_{P^\vee} (\F), \OO_{\Pic0}) , \OO_{\Pic0}\big) \Rightarrow R^{p-q}\Phi_{P^\vee}(\F).$$
In other words, to converge to $R^i\Phi_{P^\vee}(\F)$, we have to start with the $E_2$-terms given by the expressions $\mathcal{E}xt^{i +j}_{\Pic0}\big( \mathcal{E}xt^j _{\Pic0}(\R \Phi_{P^\vee} (\F), \OO_{\Pic0}) , \OO_{\Pic0}\big)$. 
But Lemma \ref{ext} implies that for all $j$ we have 
$${\rm codim}_{[\alpha]}~{\rm Supp}~\mathcal{E}xt^{i +j}_{\Pic0}\big( \mathcal{E}xt^j _{\Pic0}(\R \Phi_{P^\vee} (\F), \OO_{\Pic0}) , \OO_{\Pic0}\big) \ge i +j.$$ 
Since by hypothesis and the duality in (\ref{gd}) in order to have a non-zero expression above we must have $j \ge -k$, 
this implies that the codimension of the support of all 
$E_2$-terms is at least $i - k$. But then chasing through the spectral sequence this immediately 
implies that ${\rm codim}_{[\alpha]}~ {\rm Supp} ~ R^i\Phi_{P^\vee}(\F) \ge i-k$.
\end{proof}

\begin{remark}\label{particular}
In the projective setting, and in the general context of arbitrary 
integral transforms, the relationship between the size of the loci $V^i(\F)$ 
and the vanishing of $R^i\Phi_{P^\vee}(\R\Delta\F)$ 
is already contained in \cite{pp3}, Theorem A. The argument provided here for 
passing from (i) to (ii) is different, working also in the non-projective case and, 
together with the results of the next section, again works for general integral 
transforms. We refer to \S5 for precise statements.
\end{remark}

\section{The Generic Vanishing Index, $k$-th syzygy sheaves and the Evans-Griffith theorem}

Let $X$ be a smooth projective variety over an algebraically closed field, 
or a compact complex manifold, and let $\F$ be a coherent sheaf on $X$. We introduce the following measure of the size of  
the cohomological support loci $V^i(\F)$ for $i>0$. 

\begin{definition}\label{gv_index} 
The quantity 
$$gv (\F) : = \underset{i>0}{\rm min}~ \{ {\rm
codim}_{\Pic0}~V^i(\F)-i\}$$ 
is called the \emph{Generic Vanishing index} of $\F$.  
It is also useful
to define a local version: given any $[\alpha]\in \Pic0$, the \emph{local Generic Vanishing index} 
of $\F$ at $[\alpha]$ is 
$$gv_{[\alpha]} (\F) : = \underset{i>0}{\rm min}~ \{ {\rm
codim}_{[\alpha]}~V^i(\F)-i\}.$$   
If $V^i(\F)=\emptyset$ for all $i>0$ (respectively $[\alpha]\not\in V^i(\F)$ for all $i>0$), 
we declare $gv(\F)=\infty$ (respectively $gv_{[\alpha]}(\F)=\infty$). 
Note that by Lemma \ref{codimension}, these indices could be defined using ${\rm Supp}~R^i \Phi_{P^\vee} (\F)$ as well. Clearly $gv (\F) = {\rm inf}~ \{gv_{[\alpha]}(\F)~|~\alpha\in \Pic0\}$.
\end{definition}

The Generic Vanishing index reads a local algebraic property of the transform $\R\Phi_P(\R\Delta \F)$, namely that of being 
a syzygy sheaf (cf. the Appendix). It was already noted in \cite{pp4} that the case $m= 1$, 
i.e. torsion-freeness, characterizes the $M$-regularity condition (cf. \cite{pp1}). For $m\ge 2$, this can be interpreted 
as a characterization of \emph{higher regularity} conditions on sheaves on $X$ (or in other words of the $GV_m$-property; cf. \S5 for a more general context).

\begin{corollary}\label{syzygy_kahler}
Let $X$ be a smooth projective variety over an algebraically closed field, or a compact complex manifold, of dimension $d$, 
and let $\F$ be a coherent sheaf on $X$. Let $m \ge 0$ be an integer, or $m=\infty$, and let $[\alpha]$ be a point in $\Pic0$. The following are equivalent:

\noindent
(i) $gv_{[\alpha]}(\F) \ge  m$.

\noindent
(ii) $\widehat{\R\Delta \F}:= \R \Phi_{P} (\R\Delta \F) [d]$ is an $m$-th syzygy sheaf on $\Pic0$, in a neighborhood of $[\alpha]$.\footnote{Note 
that for $m=\infty$ this simply says that $\widehat{\R\Delta \F}$ is free in a neighborhood of $[\alpha]$.} 
\end{corollary}
\begin{proof}
By Theorem \ref{general_equiv} we know that the condition in (i) implies that $\R \Phi_{P} (\R\Delta \F) [d]$ is a sheaf
in a neighborhood of $[\alpha]$. As in its proof, we have in such a neighborhood that 
$$R^i \Phi_{P^\vee} (\F) \cong \mathcal{E}xt^i _{\Pic0}(\widehat{\R\Delta\F}, \OO_{\Pic0}).$$
Since by Lemma \ref{codimension} we have that ${\rm codim}_{[\alpha]} ~{\rm
Supp}~R^i\Phi_{P^\vee}(\F) =  {\rm codim}_{[\alpha]}~ V^i (\F)$, the result follows from Proposition \ref{higher_syzygy}.
\end{proof}   

As a consequence of the Evans-Griffith Syzygy Theorem 
(Theorem \ref{syzygy_theorem} below),  the Generic Vanishing index
of a coherent sheaf $\F$ provides a lower bound for its holomorphic
Euler characteristic in case the equivalent conditions of Corollary \ref{syzygy_kahler}
are satisfied for some $m < \infty$. This requires that there exist $i >0$ with $V^i (\F) \neq \emptyset$.

\begin{theorem}\label{chi_inequality}
Let $X$ be a smooth projective variety over an algebraically closed
field, or a compact complex manifold, and let $\F$ be a coherent 
sheaf on $X$ such that there exists $[\alpha]\in\Pic0$ with $0\le
gv_{[\alpha]}(\F)<\infty$. Then
$$\chi(\F) \ge gv_{[\alpha]}(\F) \ge  gv(\F).$$
\end{theorem}
\begin{proof} 
Note that if the equivalent
conditions of Corollary \ref{syzygy_kahler} are satisfied, then
\begin{equation}\label{rank}
{\rm rk}(\widehat{\R\Delta \F})= \chi(\F).
\end{equation} 
Indeed, by the duality formula (\ref{gd}), the dual sheaf $\mathcal{H}om_{\Pic0}(\widehat{\R\Delta \F}, \OO_{\Pic0})$ is isomorphic to $R^0\Phi_{P^\vee}(\F)$,
whose rank is the generic value of $h^0(\F\otimes\alpha^{-1})$. This
coincides with $\chi(\F)$ since, by condition (i), $V^i(\F)$ are proper 
closed subsets of $\Pic0$. Therefore from Corollary \ref{syzygy_kahler} it follows that $\widehat{\R\Delta \F}$ 
is a $gv_{[\alpha]}(\F)$-syzygy sheaf in the
neighborhood of $[\alpha]$, of rank $\chi(\F)$. By Corollary \ref{syzygy_kahler},  the hypothesis
$gv_{[\alpha]}(\F)<\infty$ means that
$\widehat{\R\Delta \F}$ is not free around $[\alpha]$.
We then apply Theorem \ref{syzygy_theorem}.
\end{proof}

\section{The holomorphic Euler characteristic of a compact K\"ahler manifold of maximal Albanese dimension}

Here we specialize the general results of the previous sections to the case when $\F$ is the canonical bundle $\omega_X$ of a compact K\"ahler manifold of maximal Albanese dimension. In this case $gv(\omega_X)\ge 0$, by the Generic Vanishing theorem of \cite{gl1}. Moreover $gv_0(\omega_X)<\infty$ (where $0$ denotes the identity point of $\Pic0$), since $0\in V^i(\omega_X)$ for all $i>0$. Hence Theorem \ref{chi_inequality} provides the 
inequality

\begin{corollary}\label{holo_chi}
Let $X$ be a compact K\"ahler manifold of maximal Albanese dimension. Then
$$\chi(\omega_X)  \ge gv_0(\omega_X)\ge gv(\omega_X).$$
\end{corollary}

\begin{proof} \emph{(of Theorem \ref{general}).}
Let $i>0$ such that ${\rm codim}~ V^i(\omega_X)-i$ computes $gv(\omega_X)$, 
and let $V$ be an irreducible component of $V^i(\omega_X)$ of maximal dimension. 
By a fundamental theorem of Green and Lazarsfeld (\cite{gl2}, Theorem 0.1), 
$V$ is a translate of a subtorus $T$ of $\Pic0$. Moreover, there
is a surjective morphism with connected fibers onto a normal compact
analytic variety $\pi:X\rightarrow Y$ such that (a) $\dim Y \le
\dim X-i$;  (b) $T\subseteq [\alpha]+\pi^* {\rm Pic}^0 (Y)$ for some $[\alpha]\in \Pic0$;  
(c) any smooth model $\tilde Y$ of $Y$ is of
maximal Albanese dimension. (Recall that the
variety $Y$ is constructed as follows: let $g:{\rm Alb}(X) \rightarrow \widehat{T}$ be 
the natural projection. Then $\pi : X \rightarrow Y$ is the Stein
factorization of the morphism $g\circ a: X \rightarrow
\widehat{T}$, where $a$ is the Albanese map of $X$.) 
Thus by (b) we have $q(\tilde Y)<\dim T$, i.e. $q(X)-q(\tilde Y)\le {\rm codim}~ V^i(\omega_X)$, while by (a) we have $\dim X - \dim Y\ge i$. In conclusion
$$q(X)-q(\tilde Y)-(\dim X - \dim \tilde Y)\le {\rm codim}~ V^i (\omega_X) - i = gv(\omega_X),$$
hence (c) and Corollary \ref{holo_chi} imply that $\pi: X \rightarrow Y$ satisfies the conclusion 
of Theorem \ref{general}.
\end{proof}

We include a variant of our main inequality, which allows for avoiding hypotheses on the ground field and the Albanese dimension in case direct computations can be made.

\begin{corollary}\label{isolated}
Let $X$ be an irregular smooth projective variety over an algebraically closed field or compact K\"ahler manifold, such that $0$ is an 
isolated point in $V^i (\omega_X)$ for all $i >0$. Then 
$$\chi (\omega_X) \ge q(X) - \dim X.$$
\end{corollary}
\begin{proof}
If the canonical bundle satisfies the Generic Vanishing condition ${\rm codim}_0~V^i (\omega_X) \ge i$ for all $i$, then  
we obtain $\chi(\omega_X) \ge gv_0(\omega_X)$ as in Corollary \ref{holo_chi}. But by assumption $V^i (\omega_X) = \{0\}$ 
in a neighborhood of $0$, hence 
$gv_0(\omega_X) = q(X) - \dim X$. Note that the hypothesis implies that the Generic Vanishing condition is equivalent to 
$q(X) \ge \dim X$. Otherwise $q(X) - \dim X <0$, but since $q(X) > 0$ and $0$ is an isolated point in all $V^i (\omega_X)$ with 
$i >0$, we have $\chi (\omega_X) \ge 0$, so the inequality still holds.
\end{proof}

\begin{remark}(Non-trivial isolated points.)\label{nontrivial} 
More generally, the same proof shows that if there exists $[\alpha]\in \Pic0$ such that $[\alpha]$ is isolated in all $V^i (\omega_X)$ with $i >0$ (meaning 
isolated in some, and not belonging to others at all, since $V^d (\omega_X) = \{0\}$), then 
$$\chi (\omega_X)\ge q(X)-\dim X + p,$$
where $p : = {\rm min} \{ k ~|~ H^k(X, \alpha^{-1}) \neq 0\}$.
This can be rephrased as follows: if $[\alpha]$ is non-trivial isolated point in $V^1(\omega_X)$, then our main inequality is improved to at least $\chi(\omega_X) \ge q(X) - \dim X +1$. Beauville has given examples of varieties of arbitrary dimension at least $2$ having such non-trivial isolated points  in $V^1 (\omega_X)$ (cf. \cite{beauville2} \S1). In the surface case, this is enough to get $p_g(X) \ge 2q(X) -2$.
\end{remark}

\section{Further results and remarks}

\noindent 
{\bf An application of Theorem \ref{general}.} 
When the holomorphic Euler characteristic is sufficiently small, Theorem \ref{general} has as an immediate consequence a strengthening of a result of Ein-Lazarsfeld \cite{el} stating that varieties of maximal Albanese dimension with $\chi(\omega_X)=0$ have the Albanese image fibered by subtori of ${\rm Alb}(X)$. We use the following notation:  let $ i_{max}$ be the maximal $i>0$ computing $gv(\omega_X)$, i.e.  such that $gv (\omega_X) = {\rm codim}~V^i (\omega_X) - i$. Then we
denote 
$$c (\omega_X): = {\rm codim}~V^{ i_{max} }(\omega_X).$$ 
As an example, one can check without difficulty that if $X=C\times D$ is the product of a curve $C$ of genus $2$ and a smooth model $D$ of a principal 
polarization in a $g$-dimensional abelian variety $A$ with $g > 2$, then 
$gv(\omega_X) = 1$ and it is computed by $i=1$ and $i=g-1$. 
Hence $i_{max}=g-1$ and $c(\omega_X)={\rm codim}~V^{ g-1 }(\omega_X)=g$, while ${\rm codim}~V^1 (\omega_X)=2$. 
Note also that $\chi (\omega_X) = 1$.

\begin{corollary}\label{ruling}
Let $X$ be a compact K\"ahler manifold of maximal Albanese dimension. If $\chi (\omega_X) \le c (\omega_X)$, then the Albanese image of $X$ is fibered by subvarieties of codimension at most $\chi (\omega_X)$ of subtori of ${\rm Alb}(X)$.
\end{corollary}

For instance, the next step beyond the Green-Lazarsfeld result is the case $\chi(\omega_X)=1$, where this says that the Albanese image of $X$ is either fibered by subtori or by divisors of subtori of ${\rm Alb}(X)$. 
For the proof, we only need to note that for a morphism
$\pi:X\rightarrow Y$ provided by Theorem \ref{general}, the dimension
of a general fiber, $\dim F_\pi = \dim X-\dim Y$, and the
dimension of the subtorus $K_\pi := \ker({\rm Alb} (X) \rightarrow {\rm
Alb}(Y))$, are subject to the inequality $\dim F_\pi \ge \dim K_\pi-\chi(\omega_X)$.
The assertion follows since the Albanese map is generically finite on $F_\pi$. 

\noindent
{\bf An example.}
Here is an example which shows that in Theorem \ref{cdf} the maximal Albanese dimension hypothesis is necessary. 
We construct a threefold $X$ with surjective Albanese map and $\chi (\omega_X) < q (X)  - 3$, such that there is no morphism with connected fibers to a smaller dimensional variety $Y$ of Albanese general type. This can be done in higher dimension as well.

Let $X = \PP^1 \times S$, where $S$ is a smooth projective surface with the following properties: $S$ is of general type, has surjective Albanese map, 
$q(S) = 2$ and $\chi (\omega_S) > 1$. Let's assume for now that such an $S$ exists and come back to this at the end. Clearly $q(X) = q(S) = 2$, and the Albanese map of $X$ is surjective.
It is immediate to check that $\chi (\omega_X) = - \chi (\omega_S)$. Thus the 
inequality $\chi (\omega_X) < q(X) - 3$ holds, since it is equivalent to $\chi(\omega_S) > 1$. 
On the other hand, there is no map $X
\rightarrow Y$ with $Y$ of Albanese general type and $\dim Y = 1,2$.
Indeed, since the Albanese map of $X$ is surjective, the induced map ${\rm Alb}(X) = {\rm Alb}(S) 
\rightarrow {\rm Alb}(Y)$ would have to surject onto the Albanese image of $Y$, which is not a subtorus, giving a contradiction.

Surfaces $S$ with the required properties can be obtained by a standard construction. 
Take any abelian surface $A$, and an ample line bundle $L$ on $A$ with $h^0 (L) \ge 2$. Then $L^{\otimes 2}$ is globally generated, so by Bertini we can consider $D \in |L^{\otimes 2}|$ a smooth irreducible divisor. Consider $\pi: S \rightarrow A$ the $2$-fold cyclic cover branched along $D$.  By the ramification formula, it is clear that $\omega_S$ is ample, so $S$ is of general type. We have $\pi_* \OO_S \cong \OO_A \oplus L^{-1}$, which combined with Kodaira Vanishing gives $q(S) = q(A) = 2$. Since $\pi $ is surjective, this in turn implies by the universal property that the Albanese map of $S$ has to be surjective as well. (In fact one can see that it is equal to $\pi$.) Finally, by relative duality we have $\pi_* \omega_S \cong \OO_A \oplus L$. This gives $\chi (\omega_S) = \chi (L) = h^0 (L) \ge 2$.

\noindent
{\bf Optimality.} 
Recall that Theorem \ref{cdf} can be rephrased as saying that 
compact K\"ahler manifolds of maximal Albanese dimension which do not admit
higher irrational pencils satisfy the inequality
\begin{equation}\label{repetition}
\chi_{hol} := \chi(\omega_X)\ge q(X)-\dim X.
\end{equation}
The natural question concerning the optimality of this lower bound
leads to intriguing problems. Equality in (\ref{repetition}) is
achieved when $\chi_{hol}=0$ (abelian varieties) and
$\chi_{hol}=1$ (desingularizations of theta divisors of principally
polarized abelian varieties), as well as for curves. 
These however seem to be the only known examples without 
higher irrational pencils. When $\chi_{hol}$ gets higher, 
the state of affairs is completely open, and there are no known examples 
even in the classical case of surfaces.\footnote{Note that by contrast the classification
of surfaces with irrational pencils and $\chi_{hol} = q(X) - 2$ is by now completely understood; cf. \cite{hp3} and \cite{pirola} for $q(X) = 3$, \cite{bnp} for $q(X) = 4$ and \cite{mlp} for $q(X) \ge 5$.
Mendes-Lopez and Pardini \cite{mlp} conjecture that there should be no examples
with $q(X) \ge 4$ and with no irrational pencils of genus at least $2$.} We propose the 
potential construction of such examples as a highly interesting problem. We suspect though 
that for $\chi_{hol} \ge 2$ equality is very rarely, if ever, achieved.

\noindent 
{\bf The setting of arbitrary integral functors.}
The notion of Generic Vanishing index and its relation with the $k$-th syzygy
property for the Fourier-Mukai transform of the Grothendieck dual is not restricted to the context of 
the Picard variety. In the 
general setting of $GV$-sheaves (cf. \cite{pp3}) with respect to 
arbitrary integral transforms between smooth varieties given by 
locally free kernels, it holds with essentially the same proofs. We will only formulate the statements
below, leaving the details to the interested reader. Extensions to objects in the derived category, 
on possibly singular varieties, will be discussed elsewhere.

Let $X$ be a compact complex manifold (respectively a smooth projective
variety over an algebraically closed field $k$) and let $Y$ be a 
complex manifold, not necessarily compact (respectively a smooth algebraic
variety over $k$). Let $\D (X)$ be the bounded derived category of
coherent sheaves on $X$, and same for $Y$ and $X\times Y$.
Any object $E$ in $\D(X\times Y)$ defines an integral transform
$$\R \Phi_E : \D (X) \rightarrow \D(Y),
~~ \R \Phi_E(\cdot) : = \R {p_Y}_* (p_X^*(\cdot) \overset{\R}{\otimes} E).$$
Let now $P$ be an object in 
$\D(X\times Y)$, and $P^\vee$ the dual object $\R{\mathcal H}om_{X\times Y}(P,\OO_{X\times Y})$.
The extension of Theorem \ref{general_equiv} under these hypotheses is:

\begin{theorem}\label{GV} 
Let $\F$ be an object in $\D(X)$ and $k\ge 0$ an integer. The following are equivalent:

\noindent
(i) ${\rm codim} ~{\rm Supp} ~R^i \Phi_{P^\vee} (\F) \ge i-k$.

\noindent
(ii) $R^i \Phi_P (\R\Delta\F) = 0$ for all $i< \dim X-k$.
\end{theorem}

In the language of \cite{pp3}, an object satisfying condition (i) is
called a \emph{$GV_{-k}$-object}. We recall that in the projective case 
this statement in contained in Theorem A of \cite{pp3}.

We now assume that $P$ is a \emph{locally free sheaf} (or more generally a perfect complex) on
$X\times Y$. One can define the Generic Vanishing index of $\F$ with respect to $P$, 
denoted $gv^P (\F)$, exactly as in Definition \ref{gv_index}. The definition of $GV_{-k}$-sheaves 
in \cite{pp3} can be extended to a full filtration of ${\rm Coh}(X)$.
\begin{definition}\label{gv_sheaves}
A coherent sheaf $\F$ on $X$ satisfies Generic Vanishing with index $k$, or simply it is a 
$GV_k$-sheaf, with respect to $P$, if $gv^P (\F) \ge k$. 
\end{definition}
\noindent
If the conditions of Theorem \ref{GV} are satisfied for
$k=0$, then the transform $\R \Phi_P (\R\Delta \F)=R^{\dim X}
\Phi_P (\R\Delta\F)[\dim X]$ is a sheaf on $Y$, denoted $\widehat{\R\Delta\F}$. 
Now, for $y\in Y$, let $P_y := P_{X\times \{y\}}$. Then $\chi(\F\otimes P_y)$ does
not depend on $y$, and will be denoted $\chi^P (\F)$. 
One has the following statements extending Theorem \ref{syzygy_kahler}.
and Theorem \ref{chi_inequality}.

\begin{theorem}
Let $0\le k\le \infty$ be an integer. Under the hypotheses above, 
$\F$ is a $GV_k$-sheaf if and only if the transform $\widehat{\R\Delta \F}$ 
is a $k$-th syzygy sheaf on $Y$. If in addition $gv^P(\F)<\infty$, then
$$\chi^P(\F)\ge gv^P(\F).$$
\end{theorem}

\section{Appendix:  $k$-th syzygy sheaves}

In this Appendix $X$ is a smooth variety over an algebraically closed field or a complex manifold.

\begin{lemma}\label{ext}
If $\F$ a coherent sheaf on $X$, then 
$$\E xt^i (\F , \OO_X) = 0 {\rm ~ for~ all~} i < {\rm codim}~{\rm Supp}~\F$$
and 
$${\rm codim}~ {\rm Supp} ~ \E xt^i (\F , \OO_X) \ge i {\rm ~ for~ all~} i.$$
\end{lemma}
\begin{proof}
This is a well-known application of the Auslander-Buchsbaum-Serre formula, cf. e.g. \cite{hl} Proposition 1.1.6(1)
and \cite{oss} Lemma II.1.1.2. 
The proof of the first assertion is given in \emph{loc. cit.} in the projective case. 
For a proof in the local case cf. e.g. \cite{bs} Proposition 1.17.
\end{proof}

\begin{definition}
A coherent sheaf $\F$ on $X$ is called a \emph{$k$-th syzygy sheaf} if
locally there exists an exact sequence
\begin{equation}\label{syzygy}
0\longrightarrow \F \longrightarrow \E_k \longrightarrow \ldots
\longrightarrow \E_1 \longrightarrow \G \longrightarrow 0
\end{equation}
with $\E_j$ locally free for all $j$. It is well known for example
that $1$-st syzygy sheaf is equivalent to torsion-free, and $2$-nd syzygy
sheaf is equivalent to reflexive. Every coherent sheaf is declared
to be a $0$-th syzygy sheaf, while a locally free sheaf is declared to
be an $\infty$-syzygy sheaf.
\end{definition}

\noindent Following \cite{eg1} and \cite{hl} \S1.1, we consider
Serre's condition for coherent sheaves.

\begin{definition}
A coherent sheaf $\F$ on $X$ satisfies property $S_k$ if for all $x$
in the support of $\F$ we have:
$${\rm depth}~ \F_x \ge {\rm min} \{k, {\rm dim}~\OO_{X, x}\}.$$
\end{definition}

The following is a combination of various standard
commutative algebra results plus a most likely well-known fact, Lemma \ref{torsion}, 
which we could not locate in the literature. 

\begin{proposition}\label{higher_syzygy}
For  a coherent sheaf $\F$ on $X$, the following are equivalent:

\noindent (a) $\F$ is a $k$-th syzygy sheaf.\\
 (b) ${\rm codim}~ {\rm Supp} ~  {\mathcal E}xt^i (\F,\OO_X) \ge i + k,
{\rm ~for~ all~} i > 0.$\\
 (c) $\F$ satisfies $S_k$.
\end{proposition}
\begin{proof}
For $k = 0$, (a) and (c) do not impose any conditions on a coherent
sheaf, while (b) also holds since, in any case
\begin{equation}\label{abs} {\rm codim}\, ~{\rm Supp}\,~{\mathcal
E}xt^i (\F,\OO_X)\ge i
\end{equation}
for any coherent sheaf $\F$, by Lemma \ref{ext}.
The equivalence of (a) and (c) is a basic result of Auslander-Bridger, 
\cite {ab} Theorem 4.25. The equivalence of (b) and (c) is
\cite{hl} Proposition 1.1.6(ii) in the case when the support of $\F$
is the entire $X$. Now for $k \ge 1$, conditions (a) and (c) clearly
imply that this is the case. We are only left with checking that (b)
also implies for $k \ge 1$ that $\F$ is supported everywhere. But
this follows from the stronger Lemma below.
\end{proof}

\begin{lemma}\label{torsion}
A coherent sheaf $\F$ on $X$ is torsion-free if and only if ${\rm
codim} ~{\rm Supp}\, ~\mathcal{E}xt^i (\F,\OO_X) > i$ for all $i>0$.
\end{lemma}
\begin{proof}
If $\F$ is torsion free then it is a subsheaf of a  locally free
sheaf $\E$. From the exact sequence $0\rightarrow \F \rightarrow
\E \rightarrow \E/\F\rightarrow 0$ it follows that, for $i >
0$, $\mathcal{E}xt^i (\F,\OO_X) \cong
\mathcal{E}xt^{i+1} (\E/\F,\OO_X)$. But then (\ref{abs}), applied to
$\E/\F$, implies that
$${\rm codim} ~{\rm Supp}\,~\mathcal{E}xt^i (\F,\OO_X) > i,
{\rm ~for~all~}i >0.$$
Conversely, since $X$ is smooth, the functor
$\R\mathcal{H}om(\,\cdot\,,\OO_X)$ is an involution. Thus
there is a spectral sequence
$$E^{ij}_2 := \mathcal{E}xt^i \Bigl( \mathcal{E}xt^j (\F, \OO_X),\OO_X\Bigr)
\Rightarrow H^{i -j} = \mathcal{H}^{i-j} \F  = \begin{cases}\F
&\hbox{if $i= j$}\cr 0&\hbox{otherwise}\cr\end{cases}.$$ 
If $ {\rm codim} ~{\rm Supp} ~\mathcal{E}xt^i (\F,\OO_X) > i$ 
for all $i > 0$, then  $\mathcal{E}xt^i \Bigl(\mathcal{E}xt^j (\F,
\OO_X),\OO_X \Bigr) = 0$ for all $i,j$ such that  $j > 0$ and $i - j
\le 0$, so the only $E^{ii}_{\infty}$ term which might be non-zero
is $E^{00}_{\infty}$. But the differentials coming into $E^{00}_p$
are always zero, so we get a sequence of inclusions
$$\F=H^0= E^{00}_{\infty} \subset\ldots \subset E^{00}_3 \subset E^{00}_2.$$
The extremes give precisely the injectivity of the natural map
 $\F \rightarrow \F^{**}$. Hence $\F$ is torsion free.
\end{proof}

\noindent
Most important for our applications is the sheaf theoretic version of the 
Syzygy Theorem of Evans-Griffith.

\begin{theorem}[\cite{eg1}, Corollary 1.7]\label{syzygy_theorem}
Let $\F$ be a $k$-th syzygy sheaf on $X$ which 
is not locally free. Then ${\rm rank}(\F) \ge k$.
\end{theorem}

\providecommand{\bysame}{\leavevmode\hbox
to3em{\hrulefill}\thinspace}

\end{document}